 \documentclass[a4paper,12pt]{amsart}
\usepackage{amsfonts}
\usepackage{amssymb}
\usepackage{ifthen}
\usepackage{graphicx}
\nonstopmode \numberwithin{equation}{section}
\usepackage[usenames]{color}
\newtheorem{thm}{Theorem}[section]
\newtheorem{lem}{Lemma}[section]
\newtheorem{cor}{Corollary}[section]
\newtheorem{prop}{Proposition}[section]

\newtheorem{cl}{Claim}[section]
\newtheorem{ca}{Case}[section]
\newtheorem{sca}{Subcase}[section]
\newtheorem{scl}{Subclaim}[section]
\newtheorem{conj}[equation]{Conjecture}

\theoremstyle{definition}
\newtheorem{defn}{Definition}[section]
\newtheorem{op}[equation]{Open Problem}
\newtheorem{ques}[equation]{Question}
\newtheorem{rem}{Remark}[section]
\newtheorem{exam}[equation]{Example}

\newcounter {own}
\def\theown {\thesection       .\arabic{own}}

\newenvironment{pf}[1][]{%
 \vskip 3mm
 \noindent
 \ifthenelse{\equal{#1}{}}%
  {{\slshape Proof. }}%
  {{\slshape #1.} }%
 }%
{\qed\bigskip}

\newcounter{alphabet}
\newcounter{tmp}

\makeatletter
\newcommand{\Ref}[1]{\@ifundefined{r@#1}{}{\setcounter{tmp}{\ref{#1}}\Alph{tmp}}}
\makeatother

\newcommand{\IR}{{\mathbb R}}

\newcommand{\diam}{{\operatorname{diam}}}

\newcommand{\dist}{{\operatorname{dist}}}

\def\be{\begin{equation}}
\def\ee{\end{equation}}

\newcommand{\ben}{\begin{enumerate}}
\newcommand{\een}{\end{enumerate}}

\newcommand{\blem}{\begin{lem}}
\newcommand{\elem}{\end{lem}}
\newcommand{\bthm}{\begin{thm}}
\newcommand{\ethm}{\end{thm}}
\newcommand{\bcor}{\begin{cor}}
\newcommand{\ecor}{\end{cor}}
\newcommand{\beg}{\begin{exam}}
\newcommand{\eeg}{\end{exam}}
\newcommand{\begs}{\begin{examples}}
\newcommand{\eegs}{\end{examples}}
\newcommand{\bdefe}{\begin{defn}}
\newcommand{\edefe}{\end{defn}}
\newcommand{\bprob}{\begin{prob}}
\newcommand{\eprob}{\end{prob}}
\newcommand{\bques}{\begin{ques}}
\newcommand{\eques}{\end{ques}}
\newcommand{\bei}{\begin{itemize}}
\newcommand{\eei}{\end{itemize}}
\newcommand{\bcon}{\begin{conj}}
\newcommand{\econ}{\end{conj}}
\newcommand{\bop}{\begin{op}}
\newcommand{\eop}{\end{op}}

\newcommand{\bca}{\begin{ca}}
\newcommand{\eca}{\end{ca}}
\newcommand{\bsca}{\begin{sca}}
\newcommand{\esca}{\end{sca}}

\newcommand{\bcl}{\begin{cl}}
\newcommand{\ecl}{\end{cl}}

\newcommand{\bscl}{\begin{scl}}
\newcommand{\escl}{\end{scl}}

\newcommand{\bcons}{\begin{conjs}}
\newcommand{\econs}{\end{conjs}}

\newcommand{\bprop}{\begin{prop}}
\newcommand{\eprop}{\end{prop}}

\newcommand{\br}{\begin{rem}}
\newcommand{\er}{\end{rem}}
\newcommand{\brs}{\begin{rems}}
\newcommand{\ers}{\end{rems}}
\newcommand{\bo}{\begin{obser}}
\newcommand{\eo}{\end{obser}}
\newcommand{\bos}{\begin{obsers}}
\newcommand{\eos}{\end{obsers}}
\newcommand{\bpf}{\begin{pf}}
\newcommand{\epf}{\end{pf}}
\newcommand{\ba}{\begin{array}}
\newcommand{\ea}{\end{array}}
\newcommand{\beq}{\begin{eqnarray}}
\newcommand{\beqq}{\begin{eqnarray*}}
\newcommand{\eeq}{\end{eqnarray}}
\newcommand{\eeqq}{\end{eqnarray*}}

\newcounter{minutes}
\divide\time by 60
\newcounter{hours}
\multiply\time by 60 \addtocounter{minutes}{-\time}

\begin{document}

\bibliographystyle{amsplain}
\title{Characterizations of John spaces}

\def\thefootnote{}
\footnotetext{ \texttt{\tiny File:~\jobname .tex,
          printed: \number\year-\number\month-\number\day,
          \thehours.\ifnum\theminutes<10{0}\fi\theminutes}
} \makeatletter\def\thefootnote{\@arabic\c@footnote}\makeatother

\author{Yaxiang  Li}
\address{Yaxiang Li,  Department of Mathematics, Hunan First Normal University, Changsha,
Hunan 410205, People¡¯s Republic of China.
} \address{  Changsha University of Science and Technology£¬Changsha 410114, Hunan, P.R.China} \email{yaxiangli@163.com}

\author[]{Matti Vuorinen}
\address{Matti Vuorinen, Department of Mathematics and Statistics, University of Turku,
FIN-20014 Turku, Finland}
\email{vuorinen@utu.fi}

\author{Qingshan Zhou${}^{\mathbf{*}}$}
\address{Qingshan Zhou, school of mathematics and big data, foshan university,  Foshan, Guangdong 528000, People's Republic
of China} \email{q476308142@qq.com}

\date{}
 \keywords{John spaces, natural condition, quasisymmetric invariance, locally quasiconvex.\\
${}^{\mathbf{*}}$ Corresponding author}

\begin{abstract}
The main purpose of this paper is to study the characterizations of John spaces. We obtain five equivalent characterizations for length John spaces. As
an application, we establish a  dimension-free quasisymmetric invariance of length John spaces.
This result is new also in the case of the Euclidean space.

\end{abstract}

\thanks{The research was partly supported by  by NNSF of
China (Nos. 11601529 and 11671127) and Hunan Provincial Key Laboratory of Mathematical Modeling and Analysis in Engineering (No.2017TP1017).}

\maketitle{} \pagestyle{myheadings} \markboth{}{}

\section{Introduction}

This work is mainly motivated by the geometric properties of length John spaces. For $a\geq 1$, a noncomplete metric space $(D,d)$ is called {\it length $a$-John with center $x_0$} if there is a distinguished point $x_0\in D$ such that for every point $x\in D$  we can find an  $a$-carrot  arc $\alpha$ joining $x$ and $x_0$. An arc is called {\it $a$-carrot} arc if for all $z\in \alpha$, we have
$$\ell(\alpha[x,z])\leq a d(z),$$
  where $d(z)=\dist(z,\partial D)$ and $\ell(\alpha[x,z])$ denotes the length of the part $\alpha$ with endpoints $x$ and $z$. The concept of length John spaces is clearly a direct generalization of the well-known John domains in Euclidean spaces, which was introduced in 1961 by F. John \cite{Jo61} in connection with his work in elasticity. The excellent references for several characterizations of John disks and John domains see \cite{RJ} and \cite{Her}.

By now the class of John domains in $\mathbb{R}^n$ and metric John domains in doubling metric spaces have been extensively studied in connection with quasiconformal analysis and Poincar\'{e} inequalities (see \cite{BKL,HaK,Heinonen,Her,HLPW,KL} and the references therein). In \cite{BKL}, Buckley et al. proved the equivalence of metric John domains and Boman domains in the abstract setting of homogeneous spaces. Moreover, Haj{\l}asz and Koskela in \cite{HaK} proved the equivalence of metric John domain with weak John property and a chain condition in a doubling metric space which shares a local  connectivity property. The main purpose of this article is to explore the equivalence conditions of John spaces which are independent of the extra metric and geometric properties. Our main result is as follows.
\begin{thm}\label{thm-1} Suppose that $D$ is a locally $(\lambda,c)$-quasiconvex, rectifiably connected, noncomplete metric space. Let $x_0\in D$. Then the following conditions are quantitatively equivalent:
\begin{enumerate}
  \item $D$ is length $a$-John with center $x_0$.
  \item   $\diam (D)\leq b d(x_0)$ and   for every $x_1\in D$, there exists a curve $\alpha$ joining $x_1$ to $x_0$ with
\be\label{li-0} \ell_k(\alpha[x_1,y])\leq b_1|\log\frac{d(y)}{d(x_1)}|+b_2\;\;\;\;\;\;\;\;\;\;\;\mbox{for all}\;\;y\in\alpha. \ee
  \item For every $x_1\in D$, we can join $x_1$ to $x_0$ by a curve $\alpha$ such that
  $$\ell_k(\alpha[x_1,y])\leq b,$$
  where either $y=x_0$ if $d(x_1)\geq \frac{1}{2}d(x_0)$ or else $y$ is the first point of $\alpha$ with $d(y)=2d(x_1)<d(x_0)$.
  \item For every $x_1\in D$, we can join $x_1$ to $x_0$ by a curve $\beta$ such that
 $\ell(\beta)\leq a|x_1-x_0|$ and $\beta$ is $a$-carrot.
 \item For any $x_1\in D$, we can join $x_1$ to $x_0$ by a curve $\alpha$ such that $\alpha$ is diameter $a$-carrot and satisfies the $\varphi$-natural condition $(\ref{li-3})$.
\end{enumerate}
The constants $a$, $b$, $b_1$, $b_2$ (not necessarily the same at each occurrence) and the function $\varphi$ depend only on each other and the constants $c$ and $\lambda$. We note that $\ell_k(\alpha)$ and $\ell(\beta)$ denote the quasihyperbolic length of $\alpha$ and the arc length of $\beta$ respectively, and their definitions will be defined in Section $2$.
\end{thm}

A careful reader might find that there is an interesting characterization for length John spaces. Indeed, (4) in Theorem \ref{thm-1} indicates that every point in a length John space can be joined to the center by a carrot arc which also satisfies a quasiconvexity condition. Evidently, it is a stronger condition but has a useful feature for length John spaces in contrast to the definition. This attractive geometric property reminds the readers of another concept which is known as {\it uniform} domains \cite{Martio-80,MS} and also that of uniform spaces which was introduced by Bonk et al. \cite{BHK}. In a sense, a length John space can be viewed as a ``one-sided" uniform space.

For the last equivalence condition, we mainly investigate the relationship between length John spaces and diameter John spaces. In $\mathbb{R}^n$, N\"{a}kki and V\"{a}is\"{a}l\"{a} \cite[Theorem $2.16$]{RJ} proved that these two kinds of domains are equivalent. On the other hand, V\"{a}is\"{a}l\"{a} \cite[Properties $3.13$ and $3.18$]{Vai2004} constructed ``a broken tube" in an infinite-dimensional separable Hilbert space, which is a diameter John domain but not length John. It is well-known that every domain in $\mathbb{R}^n$ is a natural domain, see \cite[Theorem $2.7$]{Vai-3}. Actually, this result holds in a broader setting, that is, every domain in a doubling locally quasiconvex rectifiably connected metric space with nonempty boundary is also natural. This can be shown with a similar argument as in \cite[Corollary $2.18$]{Vu} and a complete proof will be presented in our coming paper. It is easy to see that a length John space is diameter John. For the converse, we introduce a natural condition with respect to a distinguished point, for the definition see section \ref{sec-2}. We see that this natural condition together with the diameter John property are equivalent to the length John property.

As an application of our main result, we show that a length John space is invariant under the quasisymmetric mappings. This assertion is also new for John domains in the Euclidean space.

\begin{thm}\label{thm-2} Suppose $f: D \rightarrow D'$ is an $\eta$-quasisymmetric homeomorphism \; between two locally $(\lambda,c)$-quasiconvex, rectifiably connected, noncomplete metric \; spaces. If $D$ is a length $a$-John space with center $x_0$, then $D'$ is a length $a'$-John space with center $x_0'=f(x_0)$, where $a'$ depends on $\lambda$, $c$, $\eta$ and $a$.
\end{thm}

The rest of this paper is organized as follows.  In Section \ref{sec-2}, we recall some definitions and preliminary results.
Section \ref{sec-3} is devoted to the proofs of Theorems \ref{thm-1} and \ref{thm-2}.

\section{Preliminaries}\label{sec-2}

For a metric space $(D,d)$, we write $|x-y|:=d(x,y)$ for the distance between $x$ and $y$.
Throughout this paper,
balls and spheres  are written as
$$\mathbb{B}(a,r)=\{x\in D:\,|x-a|<r\},\;\;\mathbb{S}(a,r)=\{x\in D:\,|x-a|=r\}$$
and $$
\mathbb{\overline{B}}(a,r)=\mathbb{B}(a,r)\cup \mathbb{S}(a,r)= \{x\in D:\,|x-a|\leq r\}.$$

For convenience, given
spaces $(D,d)$  and $(D',d')$, a map $f:D \to D'$
and points $x$, $y$,
$z$, $\ldots$ in  $D$, we always  denote by $x'$, $y'$, $z'$, $\ldots$
the images in $D'$ of $x$, $y$, $z$, $\ldots$ under $f$,
respectively. Also, we assume that $\gamma$
denotes an arc in $D$ and $\gamma'$
 the image in $D'$ of $\gamma$
under $f$.

By a curve, we mean a continuous function $\gamma:$ $[a,b]\to D$. If $\gamma$ is an embedding of $I$,
it is also called an {\it arc}.
The image set $\gamma(I)$ of $\gamma$ is also denoted by $\gamma$. The {\it length} of $\gamma$ is denoted by
$$\ell(\gamma)=\sup\sum_{i=1}^{n}|\gamma(t_i)-\gamma(t_{i-1})|,$$
where the supremum is taken over all partitions $a=t_0<t_1<t_2\ldots<t_n=b$ of the interval $[a,b]\subset \mathbb{R}$. A metric space $(D, d)$ is called {\it rectifiably connected} if every pair of points in $D$ can be joined with a
curve $\gamma$ in $D$ with $\ell_d(\gamma)<\infty$.

Suppose $\gamma $ is a rectifiable curve or a path in a noncomplete metric space $(D,d)$, its  {\it quasihyperbolic length} is the number:

$$\ell_{k}(\gamma)=\ell_{k_D}(\gamma)=\int_{\gamma}\frac{|dz|}{d(z)}.
$$

For each pair of points $x$, $y$ in $D$, the {\it quasihyperbolic distance}
$k(x,y)$ between $x$ and $y$ is defined in the following way:
$$k(x,y)=k_D(x,y)=\inf\ell_{k}(\gamma),
$$
where the infimum is taken over all rectifiable arcs $\gamma$
joining $x$ to $y$ in $D$.

We recall the following basic estimates for quasihyperbolic distance,
first established by Gehring and Palka \cite[2.1]{GP} (see also \cite[(2.3),(2.4)]{BHK}):
\be\label{li-1} k(x,y)\geq \log\Big(1+\frac{|x-y|}
{\min\{d(x), d(y)\}}\Big)\geq \log|\frac{d(x)}{d(y)}|.\ee
In fact, more generally, we have
\be\label{li-2} \ell_{k}(\gamma)\geq \log\Big(1+\frac{\ell(\gamma)}
{\min\{d(x), d(y)\}}\Big)
\ee

We next introduce some necessary definitions.

\bdefe Let $\varphi:[0,\infty)\to [0,\infty)$ be an increasing function, but we do not require that $\varphi(0)=0$. We say that a noncomplete metric space $(D,d)$ is   {\it $\varphi$-natural with respect to $x_0$}, if there is a distinguished point $x_0\in D$ such that for every $x\in D$, one can join $x$ to $x_0$ by a curve $\alpha$ satisfying the {\it $\varphi$-natural condition}:
\be\label{li-3} \diam_k(\alpha[x,y])\leq \varphi(\frac{\diam(\alpha[x,y])}{\dist(\alpha[x,y],\partial D)})\;\;\;\;\;\;\;\;\;\mbox{for all}\;y\in\alpha,\ee where $\diam_k$ means the diameter under the quasihyperbolic metric.
\edefe

We note that each domain in the Euclidean space is $\varphi$-natural (see \cite[Theorem 2.7]{Vai-3} or \cite{Vu}). In an infinite dimensional Hilbert space, the broken tube construction  in \cite{Vai2004} provides an example of a domain, which is not natural.

\bdefe A homeomorphism $f$ from $X$ to $Y$ is said to be
$\eta$-{\it quasisymmetric} if there is a homeomorphism $\eta : [0,\infty) \to [0,\infty)$ such that
$$ |x-a|\leq t|x-b|\;\; \mbox{implies}\;\;   |f(x)-f(a)| \leq \eta(t)|f(x)-f(b)|$$
for each $t\geq 0$ and for each triple $x,$ $a$, $b$ of points in $X$.
\edefe
\bdefe A rectifiable path $\gamma$, with endpoints $x,y$, is {\it$c$-quasiconvex},
$c \geq 1$, if its length is at most $c$ times the distance between its endpoints; i.e., if $\gamma$
satisfies $$\ell(\gamma)\leq c|x-y|.$$
A metric space is {\it $c$-quasiconvex} if each pair of points can be joined by a $c$-quasiconvex
path.
Let $c \geq 1$ and $0 < \lambda \leq \frac{1}{2}$.
A noncomplete rectifiably connected metric space $(D, d)$ is said to be {\it locally $(\lambda,c)$-quasiconvex, }
if for all $x\in D$, each pair of points in $\mathbb{B}(x, \lambda d(x))$ can be joined with a $c$-quasiconvex
path.
\edefe

Let $K$ be a snowflake curve in the plane with center $(0,0)$. Let $P=(0,0,1)$, and let $G=\cup_{q\in K}[p,q]$. We see that $G$ is rectifiably connected but not locally quasiconvex.

\section{Proofs of Theorem \ref{thm-1} and Theorem \ref{thm-2}}\label{sec-3}

We give the proofs of our results by dividing this section into two subsections.
\subsection{Proof of Theorem \ref{thm-1}} We observe from the definition of length John space that the implication  $(4)\Rightarrow (1)$ is obvious. We shall prove $(1)\Rightarrow (2) \Rightarrow (3) \Rightarrow (4)$ and  $(1)\Rightarrow (5)\Rightarrow (3)$.
 Our proofs mimic those of \cite[$(3.1)$ and $(3.9)$]{GHM}, \cite[Page $330$--$331$]{Her} and \cite[Page 297]{KL}.

 $\mathbf{(1)\Rightarrow (2)}$: Pick $w\in D$ with $$|w-x_0|\geq \frac{1}{3}\diam (D).$$ Since  $D$ is length $a$-John with center $x_0$, there is a curve $\gamma$ joining $w$ to $x_0$ with $\ell(\gamma)\leq ad(x_0)$, and such that
$$\diam (D)\leq 3|w-x_0|\leq 3\ell(\gamma)\leq 3a d(x_0),$$
which implies the first assertion with $b=3a$.

Next, fix $x_1\in D$ and take a curve $\alpha$ joining $x_1$ to $x_0$ with $\ell(\alpha[x_1,y])\leq a d(y)$ for all $y\in\alpha$. We shall show that $\alpha$ satisfies inequality \eqref{li-0}. To this end we divide the proof into two cases.
\bca\label{ca-1} $\alpha\subset \mathbb{\overline{B}}(x_1,\frac{1}{2}d(x_1))$.\eca

Then for every $y\in\alpha$, we have
$$\frac{1}{2}d(x_1)\leq d(x_1)-|x_1-y|\leq d(y)\leq d(x_1)+|y-x_1|\leq \frac{3}{2}d(x_1).$$
Thus
$$\ell_k(\alpha[x_1,y])\leq 2\frac{\ell(\alpha[x_1,y])}{d(x_1)}\leq 2a\frac{d(y)}{d(x_1)}\leq 3a,$$
as desired.

\bca\label{ca-2} $\alpha$ {\rm is not contained in} $ \mathbb{\overline{B}}(x_1,\frac{1}{2}d(x_1))$.\eca

We assume that $\alpha[x_1,y]$ is not contained in $ \mathbb{\overline{B}}(x_1,\frac{1}{2}d(x_1))$ (otherwise, by Case \ref{ca-1} there is nothing to prove).
Then there is some point $w\in\alpha[x_1,y]$ with $$\ell(\alpha[x_1,w])=\frac{1}{2}d(x_1).$$ A similar argument as in Case \ref{ca-1} gives that $$\ell_k(\alpha[x_1,w])\leq 3a.$$ Moreover, we have
$$\ell_k(\alpha[w,y])= \int_{\alpha[w,y]}\frac{|dz|}{d(z)}\leq \int_{d(x_1)/2}^{\ell(\alpha[x_1,y])} \frac{ads}{s}\leq a\log\frac{2ad(y)}{d(x_1)}.$$
Therefore, we get
$$\ell_k(\alpha[x_1,y])=\ell_k(\alpha[x_1,w])+\ell_k(\alpha[w,y])\leq a\log\frac{d(y)}{d(x_1)}+(3+\log2a)a,$$
as desired.

$\mathbf{(2)\Rightarrow (3)}$: Let $x_0\in D$ with $\diam (D)\leq bd(x_0)$ for some constant $b>1$. Fix $x_1\in D$ and take a curve $\alpha$ joining $x_1$ to $x_0$ satisfying $(\ref{li-0})$.

 If $d(x_1)\geq \frac{1}{2}d(x_0)$, then appealing to $(\ref{li-0})$ with $y=x_0$, one has
$$\ell_k(\alpha[x_1,x_0])\leq b_1|\log \frac{d(x_1)}{d(x_0)}|+b_2\leq b_1\log2b+b_2,$$
the last inequality holds because $\frac{1}{2}d(x_0)\leq d(x_1)\leq \diam (D)\leq bd(x_0)$.

For the case $d(x_1)<\frac{1}{2}d(x_0)$, denote by $y$  the first point of $\alpha$ with $d(y)=2d(x_1)$, then again by $(\ref{li-0})$ one gets
$$\ell_k(\alpha[x_1,y])\leq b_1\log2+b_2.$$

$\mathbf{(3)\Rightarrow (4)}$: Fix $x_1\in D$. We prove the implication by considering three cases.

Case $\mathbf{ A}$. $|x_1-x_0|\leq \frac{\lambda}{c}d(x_0)$.

  Since $D$ is locally $(\lambda,c)$-quasiconvex, pick a curve $\beta$ joining $x_1$ to $x_0$ with
$$\ell(\beta)\leq c|x_1-x_0|\leq \lambda d(x_0),$$
for all $z\in\beta$, we have $$d(z)\geq d(x_0)-|z-x_0|\geq (1-\lambda)d(x_0),$$ which implies that
$$\ell(\beta[x_1,z])\leq \ell(\beta)\leq\frac{\lambda}{1-\lambda}d(z).$$
Hence $\beta$ is the required curve in this case.

Case $\mathbf{ B}$. $|x_1-x_0|> \frac{\lambda}{c}d(x_0)$ and $d(x_1)\geq \frac{1}{2} d(x_0)$.

It follows
 from the assumption that there is a curve $\alpha_1$ joining $x_1$ to $x_0$ satisfing $\ell_k(\alpha_1[x_1,x_0])\leq b$. Moreover, for every $w\in \alpha_1$, by $(\ref{li-1})$, we have
$$|\log \frac{d(w)}{d(x_0)}|\leq k(w,x_0)\leq \ell_k(\alpha_1[x_1,x_0])\leq b,$$ and  so $$e^{-b}d(x_0)\leq d(w)\leq e^bd(x_0).$$  We claim that $\alpha_1$ is a quasiconvex carrot arc. Indeed,  for all $w\in\alpha_1$, by $(\ref{li-2})$ one obtains that
$$\ell(\alpha_1[x_1,x_0])\leq e^{\ell_k(\alpha_1[x_1,x_0])}d(x_0)\leq e^b d(x_0)\leq e^{2b}d(w)$$ and $$\ell(\alpha_1[x_1,x_0])\leq e^b d(x_0)< \frac{ce^b}{\lambda}|x_1-x_0|,$$
as needed.

Case $\mathbf{ C}$. $|x_1-x_0|> \frac{\lambda}{c}d(x_0)$ and $d(x_1)< \frac{1}{2} d(x_0)$.

Let $n\geq 1$ be the unique integer such that $$\frac{1}{2^n}d(x_0)<d(x_1)\leq \frac{1}{2^{n-1}} d(x_0).$$
We shall construct a sequence of points $x_1,x_2,...,x_n,x_{n+1}=x_0$ and $\beta_i$ as follows. Let $\alpha_1$ be a curve  connecting $x_1$ to $x_0$ such that (3) holds and let $x_2$ be the first point of $\alpha_1$ $($when traversing $\alpha_1$ from $x_1$ towards $x_0)$ with $d(x_2)=2d(x_1)$. Denote $\beta_1=\alpha_1[x_1,x_2]$. Then $\ell_k(\beta_1)\leq b$.   Next let $\alpha_2$ be the curve joining $x_2$ to $x_0$ such that $(3)$ holds. If $d(x_2)\geq \frac{1}{2}d(x_0)$, we stop with $n=2, \beta_2=\alpha_2$ and $x_3=x_0$. Otherwise we continue the process by letting $\beta_i=\alpha_i[x_i,x_{i+1}]$ where $\alpha_i$ is the curve joining $x_i$ to $x_0$ such that $(3)$ holds, and  $x_{i+1}$ is the first point of $\alpha_i$ with $d(x_{i+1})=2d(x_i)$ . Since $d(x_i)=2^{i-1}d(x_1)$, we find that the above process stops with $i=n$, $\beta_n=\alpha_i$ and $x_{n+1}=x_0$. The above construction is essentially due to Herron \cite{Her}.

We claim that  $\beta=\bigcup_{i=1}^n \beta_i$ has the required properties.  Since $\ell_k(\beta_i)\leq b$, we first observe from the choice of $x_i$ that  for $i\in\{1,\ldots,n-1\}$, we have $$\ell(\beta_i)\leq 2bd(x_i),$$ and by $(\ref{li-2})$
$$\ell(\beta_n)\leq e^{\ell_k(\alpha_1[x_n,x_{n+1}])}d(x_0)\leq 2e^bd(x_n).$$ Let $x\in\beta$. Then there is some $1\leq j\leq n$ such that $x\in\beta_j$. Thus we get from $(\ref{li-1})$ that
$$\log \frac{d(x_j)}{d(x)}\leq k(x_j,x)\leq \ell_k(\beta_j)\leq b,\;\;\;\;\;\;\;\mbox{and so}\;\;\;\; d(x_j)\leq e^bd(x).$$
Therefore, we obtain that
$$\ell(\beta[x_1,x])\leq \sum_{i=1}^j\ell(\beta_j)\leq 2e^b\sum_{i=1}^jd(x_i)\leq 4e^bd(x_j)\leq 4e^{2b}d(x).  $$
On the other hand, since $|x_1-x_0|> \frac{\lambda}{c}d(x_0)$, we get
$$\ell(\beta)\leq 4e^{2b}d(x_0)<\frac{4c}{\lambda}e^{2b}|x_1-x_0|.$$
Hence,  the implication $\mathbf{(3)\Rightarrow (4)}$ follows.

$\mathbf{(1)\Rightarrow (5)}$: Fix $x_1\in D$. It then there follows from the preceding proof for the implication $\mathbf{(1)\Rightarrow (2)}$ that there is a length carrot curve $\alpha$ joining $x_1$ and $x_0$ satisfying $(\ref{li-0})$. Evidently, every length carrot arc is diameter carrot, it suffices to show that $\alpha$ satisfies the $\varphi$-natural condition $(\ref{li-3})$ for some increasing function $\varphi:[0,\infty)\to [0,\infty)$.

This can be seen as follows. For all $y\in \alpha$, one computes by $(\ref{li-1})$ that
\begin{eqnarray*} \diam_k(\alpha[x_1,y]) &\leq& \ell_k(\alpha[x_1,y])\leq b_1|\log \frac{d(y)}{d(x_1)}|+b_2
\\&\leq& b_1\log(1+\frac{|x_1-y|}{d(x_1)})+b_2
\\&\leq& b_1\log(1+\frac{\diam(\alpha[x_1,y])}{\dist(\alpha[x_1,y],\partial D)})+b_2.\end{eqnarray*}
So we immediately see that $\alpha$ satisfies the $\varphi$-natural condition $(\ref{li-3})$ with $\varphi(t)=b_1\log (1+t)+b_2$.

$\mathbf{(5)\Rightarrow(3)}$: Fix $x_1\in D$. Pick a diameter $a$-carrot curve $\alpha$ connecting $x_1$ to $x_0$ which satisfies the $\varphi$-natural condition $(\ref{li-3})$ for some increasing function $\varphi:[0,\infty)\to [0,\infty)$.

Then for all $y\in\alpha$, $$d(x_1)\leq d(y)+|x_1-y|\leq (1+a)d(y),$$ which implies that $$\dist(\alpha,\partial D)\geq \frac{1}{1+a}d(x_1).$$
If $d(x_1)\geq \frac{1}{2}d(x_0)$, then we have
\begin{eqnarray*}k(x_1,x_0)&\leq& \diam_k(\alpha[x_1,x_0])\leq \varphi\left(\frac{\diam(\alpha[x_1,x_0])}{\dist(\alpha,\partial D)}\right)\\&\leq& \varphi\left(a(1+a)\frac{d(x_0)}{d(x_1)}\right)\leq  \varphi(2a(1+a)):=b.\end{eqnarray*}
Pick a curve $\widetilde{\alpha}$ joining $x_1$ to $x_0$ with $\ell_k(\widetilde{\alpha}[x_1,x_0])\leq 2k(x_1,x_0)$. Thus $\ell_k(\widetilde{\alpha}[x_1,x_0])\leq 2b$ and $\widetilde{\alpha}$ is the desired curve in this case.

If $d(x_1)<\frac{1}{2}d(x_0)$, then let $y_0\in\alpha$ be the first point of $\alpha$ with $d(y_0)=2d(x_1)<d(x_0)$. By a similar argument as above for the sub-curve $\alpha[x_1,y_0]$, we obtain that $k(x_1,y_0)\leq b$. Then take a curve $\beta$ joining $x_1$ to $y_0$ with $\ell_k(\beta[x_1,y_0])\leq 2k(x_1,y_0)$. Thus $\widehat{\alpha}=\beta[x_1,y_0]\cup \alpha[y_0,x_0]$ is the required curve.

Hence the proof of this theorem is complete.\qed

\bcor Let $D\subsetneq \mathbb{R}^n$ be a length $a$-John domain with center $x_0$. Then for any $x\in D$, we can join $x$ to $x_0$ by an arc $\beta$ such that
 $\ell(\beta)\leq c|x-x_0|$ and $\beta$ is $c$-carrot with $c$ depending only on $a$.
\ecor
\subsection{Proof of Theorem \ref{thm-2}}

 Assume first that $D$ is a length $a$-John space with center $x_0$. Fix $x'\in D'$, by virtue of Theorem \ref{thm-1}, we can join $x=f^{-1}(x')$ and $x_0$ by a curve $\alpha$ such that $\alpha$ is diameter $a$-carrot and satisfies the $\varphi$-natural condition $(\ref{li-3})$ for some increasing function $\varphi:[0,\infty)\to [0,\infty)$ with $\varphi$ depending only on the hypotheses. According to \cite[Theorems $2.2$ and $2.25$]{TV}, we may assume that the extension map $f:\overline{D}\to \overline{D'}$ is also $\eta$-quasisymmetric and its inverse map $f^{-1}$ is $\eta'$-quasisymmetric with $\eta'(t)=\eta^{-1}(t^{-1})^{-1}$. To conclude that $D'$ is a length John space, again by Theorem \ref{thm-1}, we only need to show that $\alpha'=f(\alpha)$ is diameter $a'$-carrot and satisfies $\varphi'$-natural condition for some increasing function $\varphi':[0,\infty)\to [0,\infty)$ with $a'$ and $\varphi'$ depending only on the hypotheses.

We next verify the first condition in the following claim.

\bcl $\alpha'$ is a diameter $a'$-carrot arc with $a'=2\eta(a)$. \ecl
Let $z'\in\alpha'$ and $\varepsilon>0$. We choose a point $z_0'\in\partial D'$  such that $|z'-z_0'|\leq (1+\varepsilon)d(z')$. Then for every $y\in\alpha[x,z]$ we compute
\beq\label{li-4}\frac{|y'-z'|}{|z_0'-z'|}\leq \eta\left(\frac{|y-z|}{|z_0-z|}\right)\leq \eta\left(\frac{\diam(\alpha[x,z])}{d(z)}\right)\leq \eta(a).\eeq
Let $w'\in \alpha'[x',z']$ be a point such that $$\diam(\alpha'[x',z'])\leq 2|z'-w'|.$$ Then by \eqref{li-4}, we have
$$\diam(\alpha'[x',z'])\leq 2\eta(a)|z'-z_0'|\leq 2(1+\varepsilon)\eta(a)d(z').$$
Letting $\varepsilon\to 0$, we obtain that $\alpha'$ is a diameter $a'$-carrot curve   $a'=2\eta(a)$. Hence we are done.

For the second condition, we begin with two useful claims. The first one discusses the distortion of relative distance of connected sets under quasisymmetric mapping.

\bcl\label{cl-1} For every connected set $A'\subset D'$ with $\dist(A',\partial D')>0$, we have
$$\frac{\diam (A)}{\dist(A,\partial D)}\leq 6\eta'\left(\frac{\diam (A')}{\dist(A',\partial D')}\right)\;\;\;\;\;\;\;\;\;\mbox{where}\;A=f^{-1}(A').$$
\ecl
Take $x_1\in A$ and $u\in\partial D$ be points such that $|x_1-u|\leq 2\dist(A,\partial D)$. Pick another point $x_2\in A$ such that $\diam (A)\leq 3|x_1-x_2|$. Since $f^{-1}$ is $\eta'$-quasisymmetric, we compute
$$\frac{\diam (A)}{\dist(A,\partial D)}\leq  6\frac{|x_1-x_2|}{|x_1-u|}\leq   6\eta'\left(\frac{|x_1'-x_2'|}{|x_1'-u'|}\right)  \leq    6\eta'\left(\frac{\diam (A')} {\dist(A',\partial D')}\right),$$
as desired.

The next claim is related to the coarse distortion of quasihyperbolic distances under quasisymmetric mapping between two locally quasiconvex rectifiably connected noncomplete metric spaces.

\bcl\label{cl-2} There are positive numbers $c_1$ and $c_2$, depending only on the hypotheses, such that $k'(u_1',u_2')\leq c_1k(u_1,u_2)+c_2$ for all $u_1,u_2\in D$. \ecl
Thanks to \cite[Lemma $2.3$]{Vai5}, we only need to estimate $k'(u_1',u_2')$  when  $k(u_1,u_2)\leq t_0$ with $t_0$ satisfying $2\eta(e^{t_0}-1)=\frac{\lambda}{2c}$, because $(D,k)$ is a length space and evidently it is $2$-quasiconvex. Take $v_1'\in\partial D'$ with $|u_1'-v_1'|\leq 2d(u_1')$. Since $$|u_1-u_2|\leq \left(e^{k(u_1,u_2)}-1\right)d(u_1)\leq (e^{t_0}-1)d(u_1),$$   we get
$$\frac{|u_1'-u_2'|}{d(u_1')}\leq 2\frac{|u_1'-u_2'|}{|u_1'-v_1'|}\leq 2\eta\left(\frac{|u_1-u_2|}{|u_1-v_1|}\right)\leq 2\eta\left(\frac{|u_1-u_2|}{d(u_1)}\right)\leq 2\eta(e^{t_0}-1)=\frac{\lambda}{2c}.$$
Moreover, since $D'$ is locally $(\lambda,c)$-quasiconvex, there is a curve $\beta'$ joining $u_1'$ and $u_2'$ with $$\ell(\beta')\leq c|u_1'-u_2'|\leq \frac{\lambda}{2}d(u_1').$$ Thus for every $w'\in\beta'$, $$d(w')\geq d(u_1')-|u_1'-w'|\geq (1-\frac{\lambda}{2})d(u_1'),$$
which implies
$$k'(u_1',u_2')\leq \int_{\beta'}\frac{|dw'|}{d(w')}\leq \frac{\ell(\beta')}{(1-\frac{\lambda}{2})d(u_1')}\leq \frac{\lambda}{2-\lambda}.$$
This desired inequality completes the proof of this claim.

Now we are in a position to show that $\alpha'$ satisfies the $\varphi'$-natural condition \eqref{li-3}. Indeed, it follows from Claims \ref{cl-1} and \ref{cl-2} that for each $y'\in\alpha'$,
\begin{eqnarray*}\diam_{k'} (\alpha'[x',y'])&\leq& c_1\diam_k (\alpha[x,y])+c_2
\\&\leq&  c_1\varphi\left(\frac{\diam (\alpha[x,y])}{\dist(\alpha[x,y],\partial D)}\right)+c_2
\\&\leq&  c_1\varphi\left(6\eta'\left(\frac{\diam (\alpha'[x',y'])}{\dist(\alpha'[x',y'],\partial D')}\right)\right)+c_2.\end{eqnarray*}
Taking $\varphi'(t)=c_1\varphi(6\eta'(t))+c_2$, the proof of this theorem is complete.
\qed

\bcor \label{cor} Let $D\subsetneq \mathbb{R}^n$ be a length $a$-John domain with center $x_0$. If $f:D\to D'\subset \mathbb{R}^n$ is an $\eta$-quasisymmetric homeomorphism, then
$D'$ is a length $a'$-John domain with center $x'_0$ where $a'$ depends only on $a$ and $\eta$.
\ecor

\begin{rem} We remark that in \cite[Theorem 3.5]{RJ}, N\"{a}kki and V\"{a}is\"{a}l\"{a} proved that distance (or diameter) John spaces are preserved under quasisymmetric homeomorphisms, and they also showed that in $\mathbb{R}^n$ a distance (or diameter) $c$-John domain is a length $c_1$-John domain with $c_1$ depending on $c$ and $n$. These observations imply that length John domains in $\mathbb{R}^n$ are preserved under quasisymmetric homeomorphisms with the constant depending on the space dimension $n$. As  mentioned in the introduction after Theorem \ref{thm-1} that in general spaces, a distance (or diameter) John space need not be  a length  John space. So it is natural to ask whether length  John spaces are preserved under quasisymmetric homeomorphisms or not. In fact, Theorem \ref{thm-2} give an affirmative answer to this question. Moreover, Corollary \ref{cor} shows that  length John domains in $\mathbb{R}^n$ are preserved under  quasisymmetric homeomorphisms with a dimension-free control function.
\end{rem}

{\bf Acknowledgement.} The authors would like to thank the referee who has made valuable comments on this manuscript.



\begin{thebibliography}{99}



\bibitem{BHK}  {\sc M. Bonk, J. Heinonen and P. Koskela},
Uniformizing Gromov hyperbolic spaces, \textit{Ast\'{e}risque,} {\bf 270} (2001), viii+99 pp.

\bibitem{BH}  {\sc St. M. Buckley and D. A. Herron},  Uniform spaces and weak slice spaces,
\textit{Conform. Geom. Dyn.,} {\bf 11} (2007), 191--206 (electronic).

\bibitem{BHX}  {\sc S. M. Buckley, D. Herron and X. Xie}, Metric space inversions, quasihyperbolic distance, and uniform spaces, \textit{Indiana Univ. Math. J.}, {\bf 57} (2008), 837--890.

\bibitem{BKL}  {\sc S. M. Buckley, P. Koskela and G. Lu}, Boman equals John, \textit{Proc. XVI Rolf Nevanlinna colloq.}, (1996), 91--99.

\bibitem{GH}  {\sc F. W. Gehring and K. Hag}, Remarks on Uniform and quasiconformal Extension Domains, \textit{Complex Variables,} {\bf 9} (1987),
175--188.

\bibitem{GHM}  {\sc F. W. Gehring, K. Hag and O. Martio}, Quasihyperbolic geodesics in John domains, \textit{Mathematica Scandinavica,} {\bf 36} (1989), 75--92.

\bibitem{GO}  {\sc F. W. Gehring and B. G. Osgood}, Uniform domains
 and the quasi-hyperbolic metric, \textit{J. Analyse Math.,} {\bf 36} (1979),
50--74.

 \bibitem{GP}  {\sc F. W. Gehring and B. P. Palka}, Quasiconformally
homogeneous domains, \textit{J. Analyse Math.,} {\bf 30} (1976),
172--199.



\bibitem{GK}  {\sc C.Y. Guo, and P. Koskela},  Generalized John disks, \textit{ Cent. Eur. J. Math.}, {\bf 12} (2014),  349--361.

\bibitem{GKT}  {\sc C.Y. Guo, P. Koskela and J. Takkinen}, Generalized quasidisks and conformality, \textit{Publ. Mat.}, {\bf 58} (2014),  193--212.


\bibitem{HaK}  {\sc P. Haj{\l}asz and P. Koskela}, Sobolev met Poincar\'{e}, \textit{Mem. Amer. Math. Soc.,} {\bf 145} (2000),
 688, x+101 pp.
\bibitem{Heinonen} {\sc J. Heinonen}, Quasiconformal mappings onto John domains, Rev. Mat. Iberoamericana, {\bf 5} (1989), 97--123.

\bibitem{Hei}  {\sc J. Heinonen},
Lectures on analysis on metric spaces,
\textit{Springer-Verlag, Berlin-Heidelberg-New York,} 2001.

\bibitem{Her}  {\sc D. Herron}, John domains and the quasihyperbolic metric, \textit{Complex Variables and Elliptic Equations}, {\bf 39} (1999),  327--334.

   \bibitem{HLPW} {\sc M. Huang, Y. Li, S. Ponnusamy and X. Wang}, The quasiconformal subinvariance property of John domians in $\IR^n$ and its applications,
\textit{Math. Ann.,} {\bf 363} (2015), 549--615.

\bibitem{Jo61} {\sc F. John}, Rotation and strain,  \textit{Comm. Pure. Appl. Math.}, {\bf 14} (1961), 391--413.

\bibitem{KL} {\sc K. Kim and N. Langmeyer}, Harmonic measure and hyperbolic distance in John disks,  \textit{Math. Scand}, {\bf 83} (1998), 283--299.


\bibitem{Martio-80}  {\sc O. Martio},
Definitions of uniform domains, \textit{Ann. Acad. Sci. Fenn. Ser. A
I Math.,} {\bf 5} (1980), 197--205.

\bibitem{MS}  {\sc O. Martio and J. Sarvas},
Injectivity theorems in plane and space, \textit{Ann. Acad. Sci.
Fenn. Ser. A I Math.,} {\bf 4} (1978), 383--401.

\bibitem{RJ}  {\sc R. N\"{a}kki and J. V\"{a}is\"{a}l\"{a}}, John disks,
\textit{Expo. Math.} {\bf 9} (1991), 3--43.

\bibitem{TV}  {\sc P. Tukia and J. V\"{a}is\"{a}l\"{a}}, Quasisymmetric embeddings of metric spaces,
\textit{Ann. Acad. Sci. Fenn. Ser. A I Math.,} {\bf 5} (1980),
97--114.

\bibitem{Vai-3}  {\sc J. V\"{a}is\"{a}l\"{a}}, Free quasiconformality
in Banach spaces. III, \textit{Ann. Acad. Sci. Fenn. Ser. A I Math.,}
{\bf 17} (1992), 393--408.

\bibitem{Vai5} {\sc J. V\"{a}is\"{a}l\"{a}}, The free quasiworld,
Freely quasiconformal and related maps in Banach spaces,
\textit{Quasiconformal geometry and dynamics $($Lublin 1996$)$,
Banach Center Publications, Vol. 48, Polish Academy of Science,
Warsaw, ed. by Bogdan Bojarski, Julian  {\L}awrynowicz, Olli Martio, Matti Vuorinen and
J\'ozef Zaj\c{a}c},
1999, 55--118.

\bibitem{Vai2004} {\sc J. V\"{a}is\"{a}l\"{a}},Broken tubes in Hilbert spaces, \textit{Analysis}, {\bf24}, (2004), 227--238.


\bibitem{Vu} {\sc M. Vuorinen}, Capacity densities and angular limits of quasiregular mappings. \textit{ Trans.
Amer. Math. Soc.}, {\bf 263}, (1981), 343--354.


\end{thebibliography}
\end{document}